\documentclass[11pt,oneside]{amsart}

\usepackage{graphicx, amsmath, amssymb, amscd, mathtools, hyperref, amsthm, euscript, 
amsfonts,bm,color}  
\DeclareMathAlphabet{\mathpzc}{OT1}{pzc}{m}{it}

\usepackage{amsfonts,amssymb,amsmath,amscd}

\usepackage{wrapfig}
\usepackage{psfrag, graphicx, pgf, float}

\newcommand{\br}{\mathbb R}

\newcommand{\G}{\Gamma}

\newcommand{\be}{\begin{enumerate}}
\newcommand{\ee}{\end{enumerate}}
\newcommand{\bi}{\begin{itemize}}
\newcommand{\ei}{\end{itemize}}
\newcommand{\bpm}{\begin{pmatrix}}
\newcommand{\epm}{\end{pmatrix}}
\newcommand{\ba}{\backslash}
\newcommand{\op}{\operatorname}

\newcommand{\PSL}{\op{PSL}}
\newcommand{\z}{\mathbb{Z}}
\newcommand{\ga}{\gamma}
\newcommand{\q}{\mathbb Q}
\newcommand{\Ga}{\Gamma}
\newcommand{\T}{\op{T}}
\newcommand{\bH}{\mathbb{H}}
\newtheorem{Thm}{Theorem}
\newtheorem{thm}{Theorem}
\newtheorem{Def}[Thm]{Definition}
\newtheorem{Q}[Thm]{Question}

\title{Euclidean traveller in hyperbolic worlds}

\author{Hee Oh}

\begin{document}

\maketitle
We will discuss all possible closures of a Euclidean line in various geometric spaces.
Imagine the Euclidean traveller, who travels only along a Euclidean line. She will be travelling to many different geometric worlds, and
our question will be
$$\text{{\it what places does she get to see in each world?}}$$

\begin{figure}[ht]\begin{center}
  \includegraphics [height=4cm]{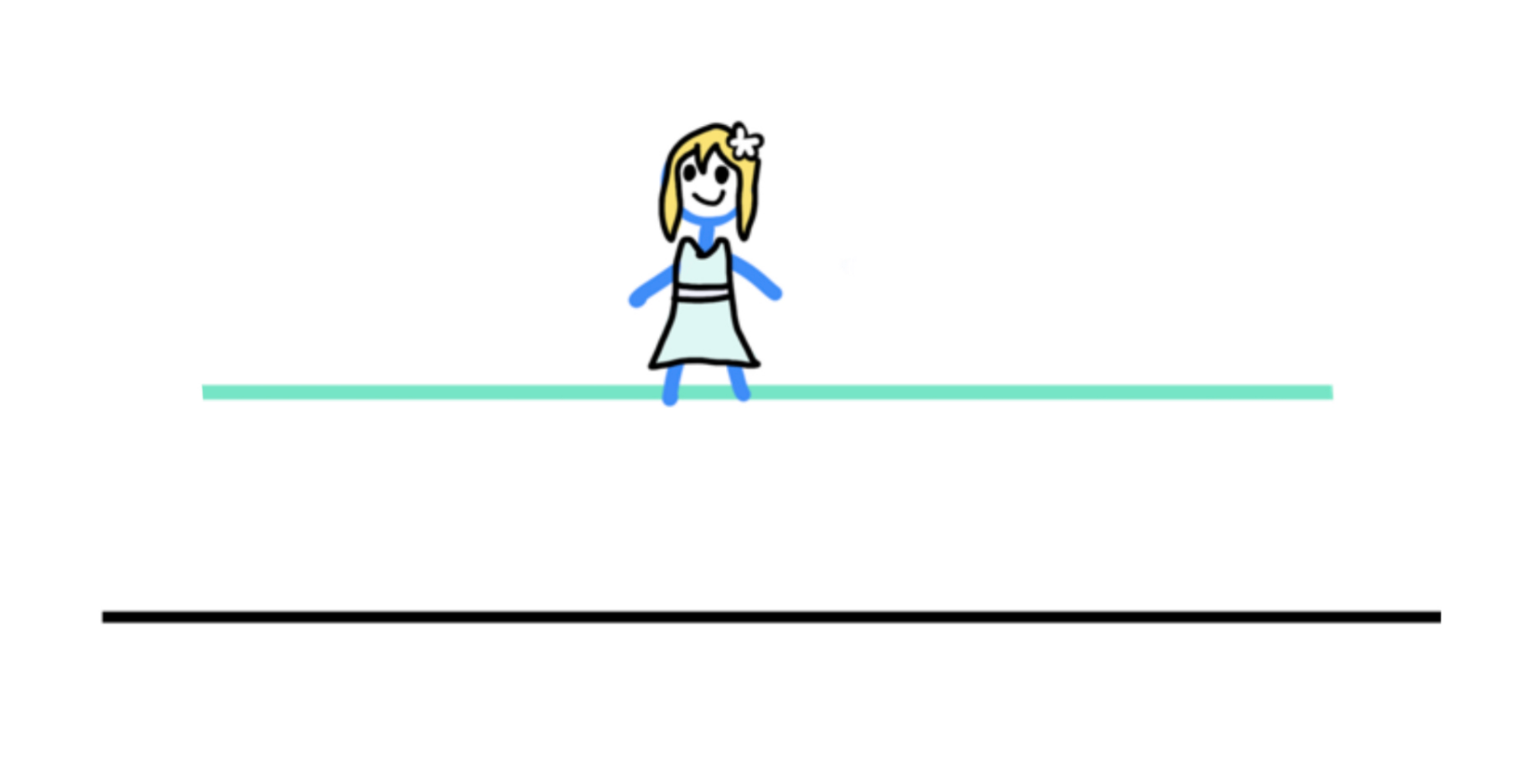}\end{center}
  \caption{A Euclidean traveller}
\end{figure}

Here is the itinerary of our Euclidean traveller:
\begin{itemize}
    \item In 1884, she  travels to the torus of dimension $n\ge 2$, guided by Kronecker. 
    \item In 1936, she travels to the world, called a
 closed hyperbolic surface, guided by Hedlund \cite{He}.
\item In 1991, she then travels to a closed hyperbolic manifold of higher dimension $n\ge 3$ guided by Ratner \cite{Ra}. 
\item Finally, she adventures into hyperbolic manifolds of infinite volume guided by Dal'bo \cite{Da} in dimension $2$ in 2000, by McMullen-Mohammadi-O. \cite{MMO2} in dimension $3$ in 2016 and by Lee-O. \cite{LO} in all higher dimensions in 2019.
\end{itemize}

\section*{Rotations of the circle}
As a warm up, she will first do her exercise of jumping on the circle $\mathbb S^1$,
which may be considered as the $1$-dimensional torus $\mathbb T^1$.

The circle $\mathbb S^1$ can be presented as
$$\mathbb S^1=\{z\in \mathbb C: |z|=1\} =\{z=e^{2\pi i x}:x\in \br\}.$$
In additive notation, 
$$\mathbb S^1= \mathbb Z\ba \br .$$
These two models are isomorphic to each other by the logarithm map
$z=e^{2\pi i x}\mapsto x\;(\text{ mod }1)$.

\begin{tabular}{cl}  
         \begin{tabular}{c}
           \includegraphics[height=3cm]{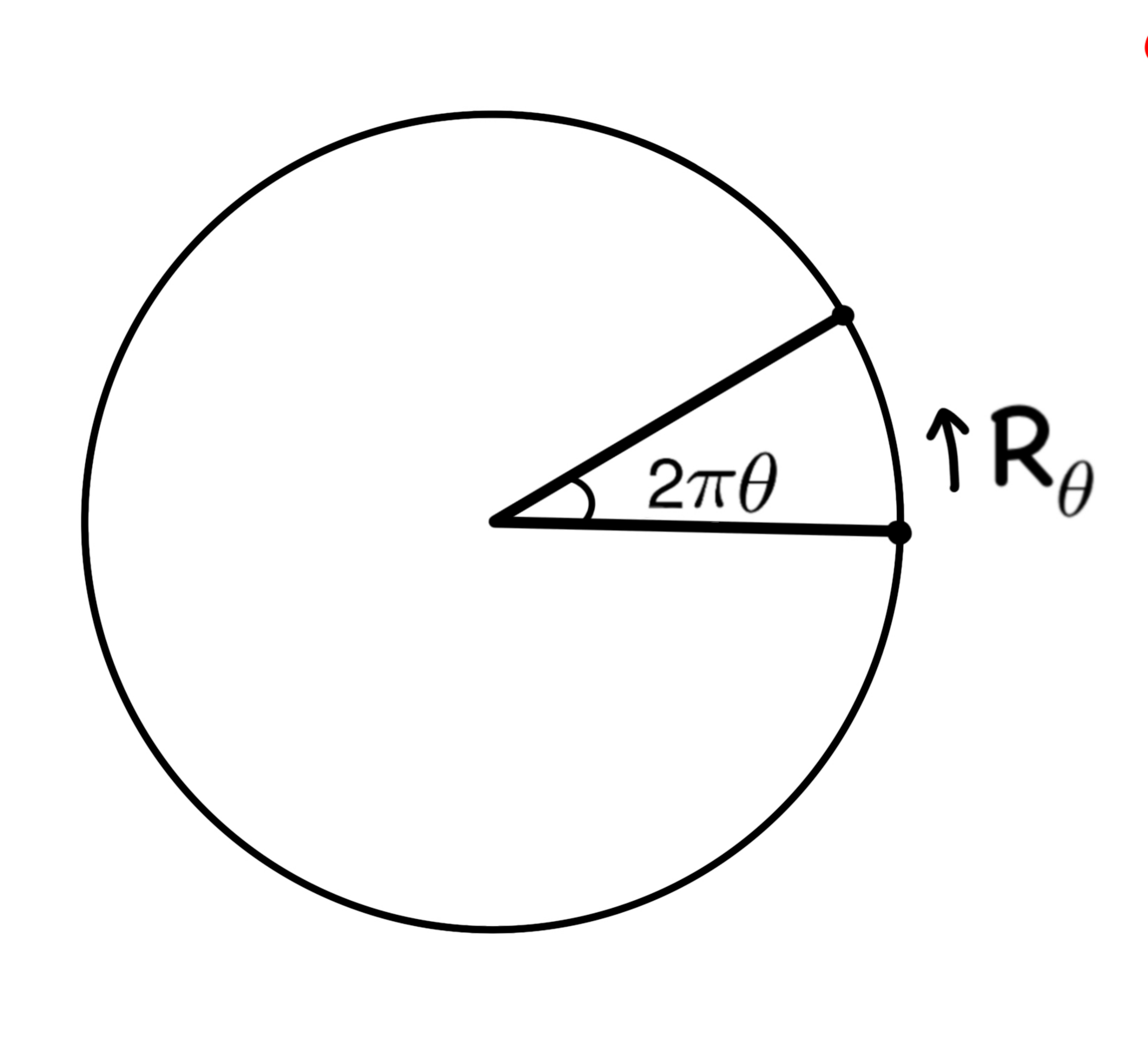}
           \end{tabular}
           & \begin{tabular}{r}
             \parbox{0.4\linewidth}{Let $\mathsf{R}_\theta$ denote the rotation of $\mathbb S^1$ by angle
$2\pi \theta$. } 
         \end{tabular} \\
\end{tabular}

The map $\mathsf R_\theta$ is  respectively given  by
$$ \mathsf{R}_\theta(z)= ze^{2\pi i \theta} \text{ and } \mathsf{R}_\theta(x)= x+\theta \quad\text{(mod $\z$)}$$
in multiplicative and additive models of $\mathbb S^1$.
The orbit of $z=e^{2\pi i x}\in \mathbb S^1$ under iterations of $\mathsf{R}_\theta$ is respectively equal to
 $$ \{ze^{2\pi i n \theta}: n\in \mathbb Z\}\text{ and }
 \{{{} x+n\theta} \;\; (\text{mod } 1): n\in \mathbb Z\}.$$
\begin{thm}\label{Ro1}
Let $\theta\in \br$. Any orbit of $\mathsf{R}_\theta$ is 
closed or dense in $\mathbb S^1$, depending on whether
$ \theta$ is rational or not.
\end{thm}

If our traveller keeps jumping by an irrational distance $\theta$, she is guaranteed to see all the places in the circle $\mathbb Z\ba \br $.

\section*{Euclidean lines on the torus}
\subsection*{$2$-torus} She travels to the $2$-torus 
$\mathbb T^2= \z^2\ba \br^2$. 
Let $$\pi: \br^2\to \mathbb T^2= \z^2\ba \br^2$$ denote the canonical quotient map. A Euclidean line in $\mathbb T^2$ is the image of a line in $\br^2$ under $\pi$. For a non-zero vector $(\omega_1, \omega_2)\in \br^2$, we denote by  
$$ L_{\omega_1,\omega_2}=\pi (\br(\omega_1, \omega_2))$$
the image of the unique line passing through $(\omega_1, \omega_2)$ and the origin $(0,0)$.
The slope of the line $\br (\omega_1, \omega_2)$ is equal to $\theta={\omega_2}/{\omega_1}$.

 What is the closure of the line $L_{\omega_1, \omega_2}$ in $\mathbb T^2$, or in other words, what places does our Euclidean traveller get to see in $\mathbb T^2$ if she travels along the line $L_{\omega_1, \omega_2}$?

  \begin{figure}[ht] \begin{center}
  \includegraphics [height=4cm]{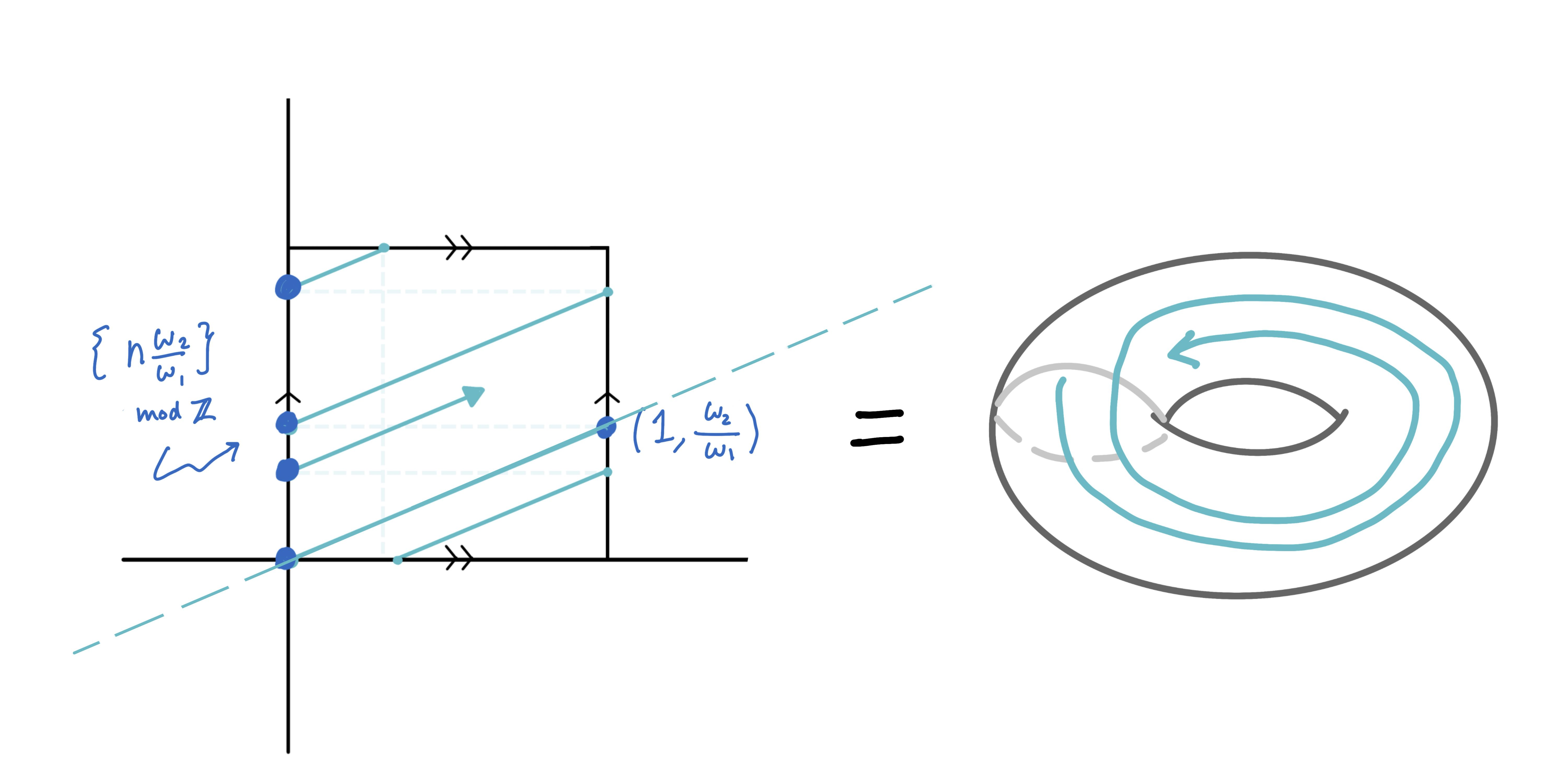} \end{center} 
  \caption{Torus} \label{T2} \end{figure}

 It is useful to consider the unit square 
 $I^2=[0,1]\times [0,1]$ which is a fundamental domain of $\mathbb Z^2\ba \br^2$. When we identify the pairs of opposite sides labelled by
 $a$ and $b$ in Figure \ref{T2} using
 the side pairing transformations $\ga_a=(1,0)$ and $\ga_b=(0,1)$ respectively, we get the torus $\mathbb T^2= \Gamma\ba \br^2$ where $\Gamma=\mathbb Z^2$ is the subgroup generated by $\ga_a$ and $\ga_b$. The distribution of the line $L_{\omega_1, \omega_2}$ can be understood by examining the line $\br({\omega_1, \omega_2})$ inside $I^2$ modulo the action of $\mathbb Z^2$. In Figure \ref{T2},
 when our traveller, walking on the blue line $\br_{\omega_1, \omega_2}$ with slope $0<\theta<1$, reaches the boundary of $I^2$ at $(1, \theta)$, she gets instantly jumped to $(0, \theta)$ by the transformation $\ga_a^{-1}=(-1,0)$, which also moves the line $\br(\omega_1, \omega_2)$ to the line  with the same slope but passing through $(0, \theta)$. She then continues to walk on this new
line in $I^2$ until she reaches the boundary of $I^2$.
We can observe that the places on the circle $\{0\}\times [0,1)=\{0\}\times \z\ba \br$ that she is visiting
are precisely given by the orbit $\{n\theta\; (\text{ mod } 1 ): n\in \mathbb N\}$ of the rotation $\mathsf{R}_{\theta}$.
Therefore by Theorem \ref{Ro1}, if the slope $\theta=\omega_2/\omega_1$ is a rational number,
she visits only finitely many points in $\{0\}\times \z\ba \br $, which means the line $L_{\omega_1, \omega_2}$ is periodic and hence closed in $\mathbb T^2$. Otherwise, the places she visits in 
  $\{0\}\times \z\ba \br $ are dense, which means that the line is dense in $\mathbb T^2$.

  \begin{thm}[Kronecker, 1884]
 Any  Euclidean line in $\mathbb T^2$ 
is  closed or dense, depending on whether its slope is rational or not.
\end{thm}

\subsection*{$n$-torus}
For any $n\ge 2$, the $n$-torus  $\mathbb T^n$ can be presented as
$\z^n\ba \br^n$, and a line in $\mathbb T^n$ is the image of a line under
the quotient map $$\pi: \br^n\to \mathbb T^n=\z^n\ba \br^n.$$
For $(\omega_1, \cdots, \omega_n)\in \br^n$,
let $L_{\omega_1, \cdots ,\omega_n}$ denote the line in $\mathbb T^n$ which is the image of the Euclidean line $\br(\omega_1, \cdots , \omega_n )$.

 \begin{figure}[ht]  \begin{center}
  \includegraphics [height=4cm]{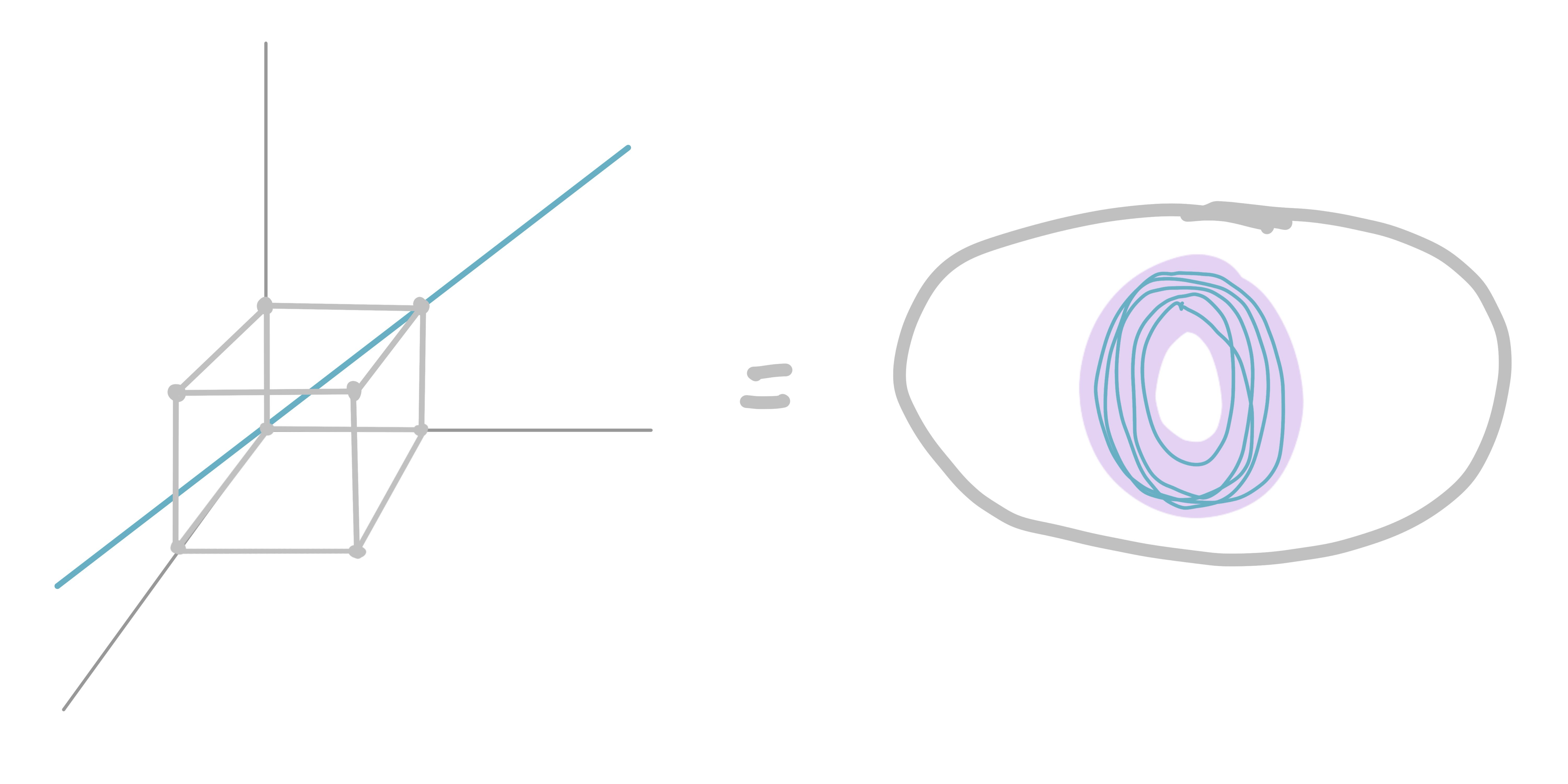} \end{center}
  \caption{Line in $\mathbb T^n$}\label{TT}
\end{figure}

    Unlike the dimension two case, the closed or dense dichotomy for a line is not true
    any more in $\mathbb T^n$, $n\ge 3$, as $\mathbb T^n$ contains many lower dimensional tori. In Figure \ref{TT}, the  blue line $L$ is contained in a two dimensional linear 
    subspace $V<\br^3$ such that
    $\pi(V)$ is closed in $\mathbb T^2$ and $\overline{L}=\pi(V)=(V\cap \z^3)\ba V$. 
    However, the lower dimensional tori are the only other possibilities (cf. \cite{KB}):

\begin{thm}[Kronecker] \label{kr1} For any non-zero $(\omega_1, \cdots, \omega_n)\in \br^n$,
the closure of the line $L_{\omega_1, \cdots, \omega_n}$ is a $k$-dimensional  subtorus of $\mathbb T^n$ where $ k=\op{dim}_\q  \sum_{i=1}^n \q \omega_i$. That is, there exists a $k$-dimensional linear subspace $V<\br^n$ such that
$$\overline{L_{\omega_1, \cdots, \omega_n}}=\pi(V)=(V\cap \z^n)\ba V.$$
\end{thm}

A general Euclidean torus is defined as the quotient 
$$\Gamma\ba \br^n$$
for some discrete cocompact subgroup $\Gamma $ of $\br^n$; a discrete subgroup $\G$ of $\br^n$ is called cocompact if the quotient space $\Gamma\ba \br^n$ is compact.
It is not hard to prove that any discrete cocompact subgroup of $\br^n$ is of the form
 $$\Gamma=\sum_{i=1}^n \z v_i$$
 for some basis $v_1, \cdots, v_n$ of $\br^n$.

\begin{thm}[Kronecker] \label{Kro}
For any line $L\subset \Gamma\ba \br^n$ (not necessarily passing through the origin), there exists a linear subspace $V<\br^n$ such that
$$\overline{L}= (V\cap \Gamma)\ba V \quad \text{up to translations}. $$
\end{thm}
This seemingly more general theorem follows easily from Theorem \ref{kr1} using the fact that $\overline{L +v}=\overline{L} +v$ for any $v\in \z^n\ba \br^n$ and
$\Gamma = \z^n g$ where $g$ is an element of $\op{GL}_n(\br)$ whose row vectors are given by $v_1, \cdots, v_n$.

\section*{Closed hyperbolic surfaces}
Our Euclidean traveller now wants to explore a world called 
a closed hyperbolic surface. A closed hyperbolic surface will be defined as a {\it quotient} of the hyperbolic plane $\bH^2$.
\subsection*{Hyperbolic plane} 
The hyperbolic plane $\bH^2$ is the unique simply connected two dimensional  complete Riemannian manifold  of constant sectional curvature $-1$. Instead of this fancy description, we will
 be using a very explicit model, called
 the Poincar\'e upper half plane model of $\bH^2$. That is,
$$\bH^2=\{(x,y)\in \br^2: y>0\}$$
with the hyperbolic metric given by $ds=\tfrac{\sqrt{dx^2+dy^2}}{y}$.
This means that 
the hyperbolic distance between $p,q\in \mathbb H^2$ is defined as 
 $$d(p, q)=\inf \left\{\int_0^1 \tfrac{\|c'(t)\| }{y(t)} dt\right\} $$ where $c:[0,1]\to \bH^2$, $c(t)=(x(t), y(t))$, ranges over all differentiable curves with $c(0)=p$ and $c(1)=q$.
Because the hyperbolic distance is the Euclidean distance scaled by the Euclidean height of the $y$-coordinate,  the hyperbolic distance between two points in $\bH^2$ is  larger
(resp. smaller) than their Euclidean distance if their Euclidean height is small
(resp. large).

\begin{figure}[ht] \begin{center}
 \includegraphics [height=3cm]{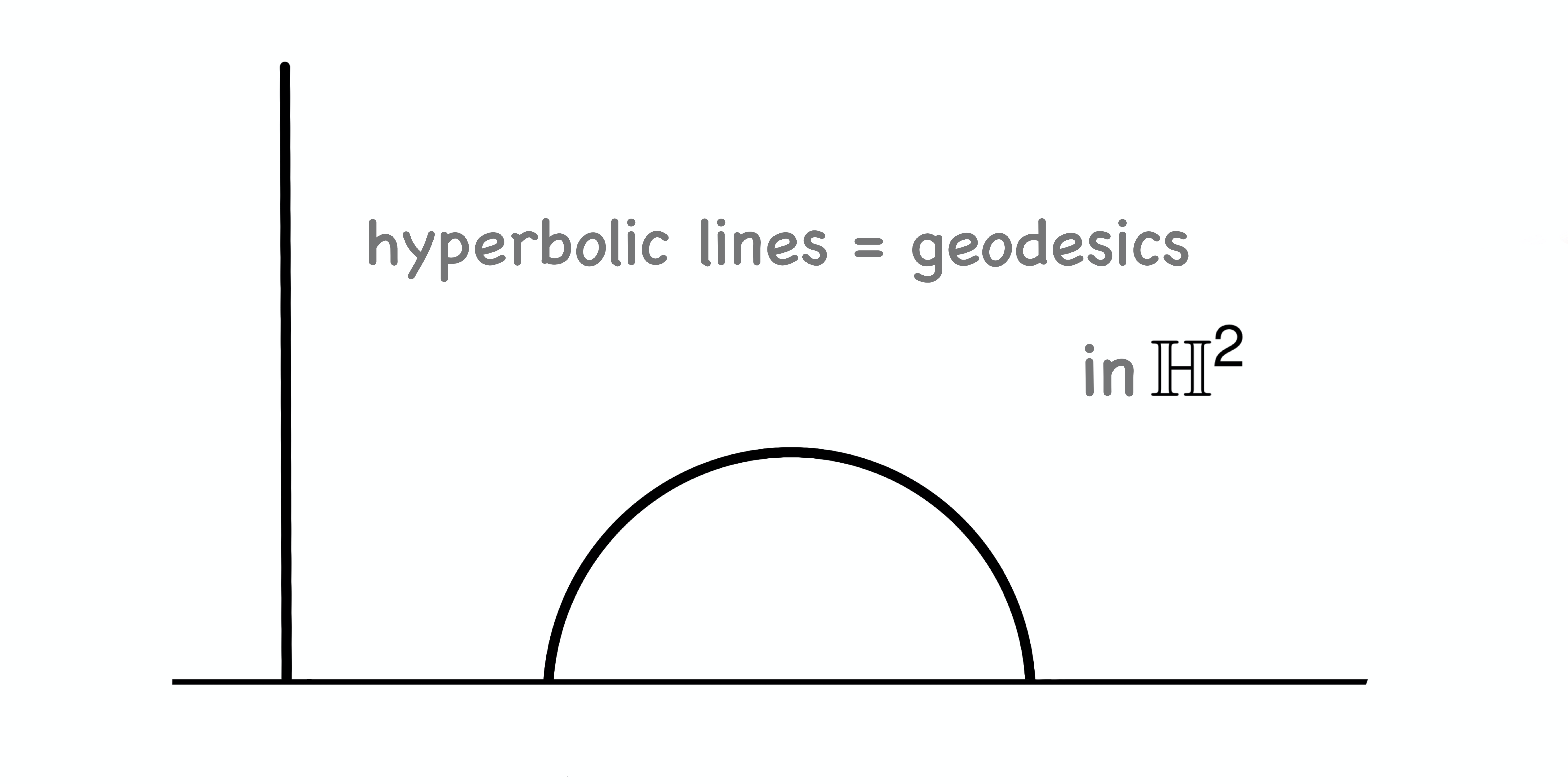}\end{center}
 \caption{Upper half-plane model of $\bH^2$}
\end{figure}

 Geodesics, that is, distance minimizing curves, in this upper half plane model are half-lines or semicircles, perpendicular to the $x$-axis.
In other words, to travel from a point $p$ to $q$ in the fastest way, one has to take the route given by the unique circle (or a line) passing through $p$ and $q$, perpendicular to the $x$-axis. 

Another useful model is the Poincar\'e unit disk model in which 
$\bH^2=\{x^2+y^2<1\}$ is the open unit disk in $\mathbb R^2$ and the hyperbolic metric is given by
$ds=\tfrac{2\sqrt{dx^2+dy^2}}{1-(x^2+y^2)}$.

\begin{figure}[ht] \begin{center}
 \includegraphics [height=3cm]{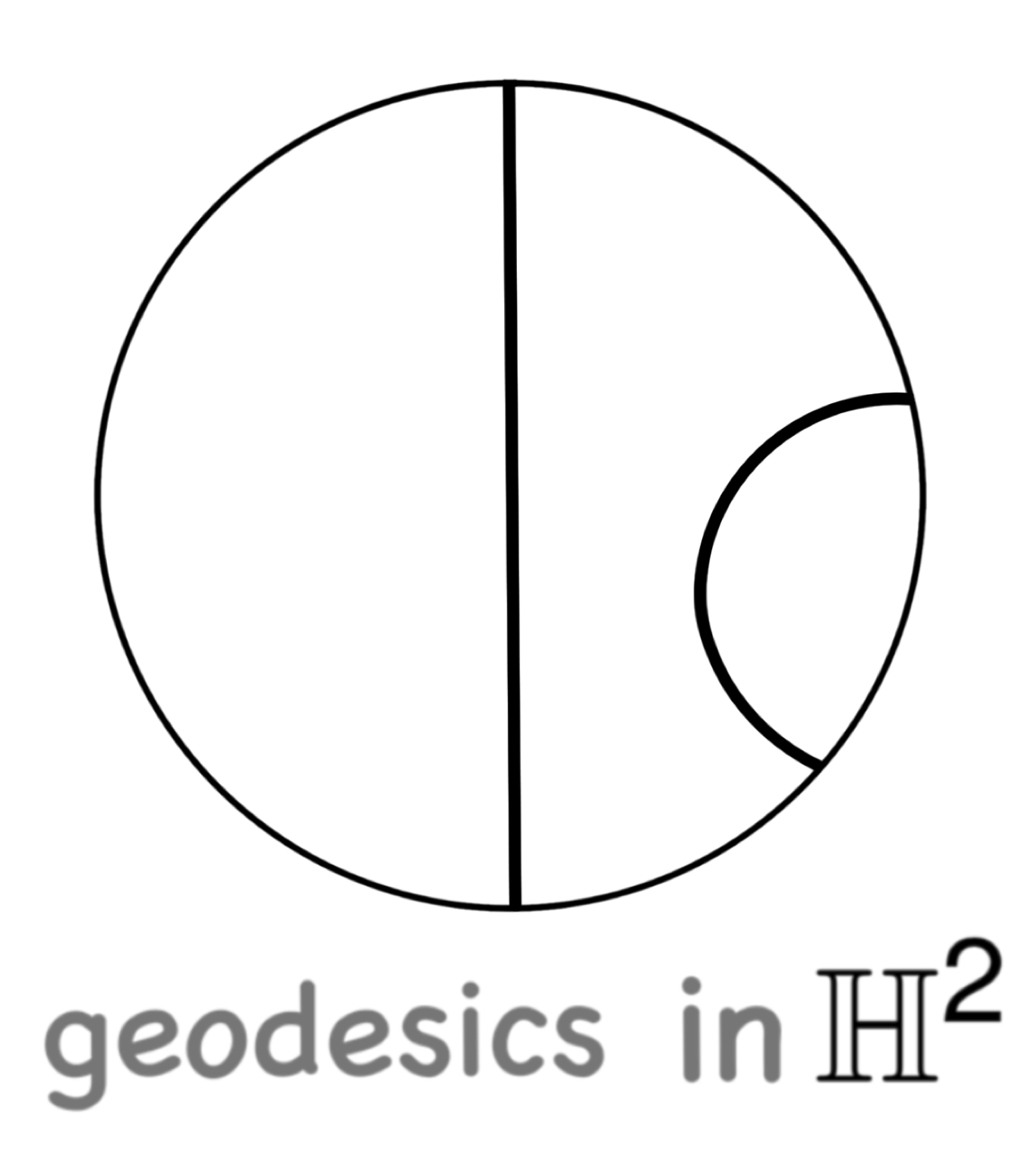}\end{center}
 \caption{Upper half-plane model of $\bH^2$}
\end{figure}
{Geodesics in this model are 
lines or semicircles meeting the boundary $\mathbb S^1$ orthogonally.}


    Recall that
  \bi\item  a Euclidean $2$-torus is given by a quotient
    $$\Gamma\ba \br^2$$ where $\Ga$ is a  discrete cocompact subgroup of $\br^2$.
   \ei 
   
 In analogy, we {\it wish} to be able to say that
    \bi 
\item   a closed hyperbolic surface
    is  given by a quotient $$\Gamma\ba \bH^2$$
where   $\Gamma$ is discrete cocompact subgroup of $\bH^2$.
    \ei 
Alas, the hyperbolic plane $\bH^2$ is not a group which makes this statement nonsense.
However, $\bH^2$ is ``almost'' the same as the group  $\op{Isom}^+(\bH^2)$ of orientation preserving isometries of $\bH^2$.

\subsection*{Isometry group of $\mathbb H^2$}
On the upper half-plane model $\bH^2=\{z=x+iy: y>0\}$, now considered as the set of complex numbers with positive imaginary parts,
the group $\op{PSL}_2(\br)$ acts by linear fractional transformations: 
$$\begin{pmatrix} a&b\\ c& d \end{pmatrix} z=\frac{az+b}{cz+d};$$
since $\op{Im} \frac{az+b}{cz+d}=\frac{\op{Im}(z)}{|cz+d|^2}$,
this action preserves $\bH^2$. We can also check that this action preserves the hyperbolic metric of $\bH^2$. Therefore every element of $\PSL_2(\br)$
is an isometry of $\bH^2$. Moreover, it turns out that every orientation preserving isometry arises in this manner, yielding the identification
 $$\op{PSL_2}(\br) =\op{Isom}^+(\bH^2) .$$

It is easy to see that this action of $\op{PSL_2(\br)}$ on $\bH^2$ is transitive with the stabilizer of $i$ being equal to the rotation group $\op{SO}(2)$.
Therefore the orbit map $\PSL_2(\br)\to \bH^2$, $g\mapsto g(i)$, induces the identification
$$\op{PSL_2}(\br)/\op{SO}(2)=\bH^2.$$

So modulo the compact subgroup $\op{SO}(2)$, the hyperbolic plane $\bH^2$ is equal to its isometry group $\PSL_2(\br)$. 

\subsection*{Closed hyperbolic surfaces}
\begin{Def} 
A closed hyperbolic surface is a quotient
$$\Gamma\ba \bH^2$$ where $\Gamma$ is a discrete (torsion-free) cocompact  subgroup of $\op{PSL}_2(\br)$. 
\end{Def} By the quotient $\Gamma\ba \bH^2$, we mean the set of equivalence classes of elements of $z\in \bH^2$ where $[z_1]=[z_2]$ if and only if $z_1=\gamma z_2$ for some $\gamma\in \Gamma$. The discreteness of $\Gamma$
implies that $\Gamma\ba \bH^2$ is locally $\bH^2$ and
the cocompactness of $\Gamma$ in $\PSL_2(\br)$ implies that $\Gamma\ba \bH^2=\Gamma\ba \PSL_2(\br)/\op{SO}(2)$ is compact.

 \begin{figure}[ht] \begin{center}
  \includegraphics [height=4cm]{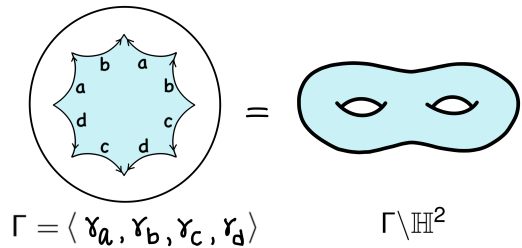}
 \caption{Hyperbolic octagon} 
  \end{center} \label{oct} 
 \end{figure}
Does there exist a closed hyperbolic surface? Equivalently, does there exist a discrete cocompact subgroup of $\op{PSL}_2(\br)$?  To present an example,
consider the  hyperbolic regular octagon $\mathsf O$ as illustrated in Figure \ref{oct}, where each side is a hyperbolic geodesic segment of same length and angles between them are $45^\circ$.

For the two sides of the octagon labelled by $a$, there exists an isometry
$\ga_a\in \PSL_2(\br)$ which moves one to the other. Similarly, we have
$\ga_b, \ga_c,\ga_d\in \PSL_2(\br)$ for labels $b, c, d$ respectively. Let
$$\Gamma=\langle \ga_a, \ga_b, \ga_c, \ga_d\rangle <\PSL_2(\br) $$ be the subgroup generated by these four side-pairing transformations.
Then the hyperbolic octagon $\mathsf O$ is a
fundamental domain for the action of $\Gamma$ in $\bH^2$, which implies that
$\Gamma$ is a discrete cocompact subgroup of $\PSL_2(\br)$.
The closed hyperbolic surface $\Gamma\ba \bH^2$ is  what we obtain by gluing the four pairs of edges of the hyperbolic octagon according to labels; it is topologically a two-holed torus, or genus-two surface.

Any closed hyperbolic surface is topologically a closed surface with genus at least $2$. Conversely, the {\it Uniformization} theorem says that
any closed surface with genus at least $2$ can be realized as a closed hyperbolic surface.
Indeed, for any $g\ge 2$, the space of all marked
closed hyperbolic surfaces of genus is homeomorphic to $\br^{6g-6}$.

 \begin{figure}[ht] \begin{center}
  \includegraphics [height=3.5cm]{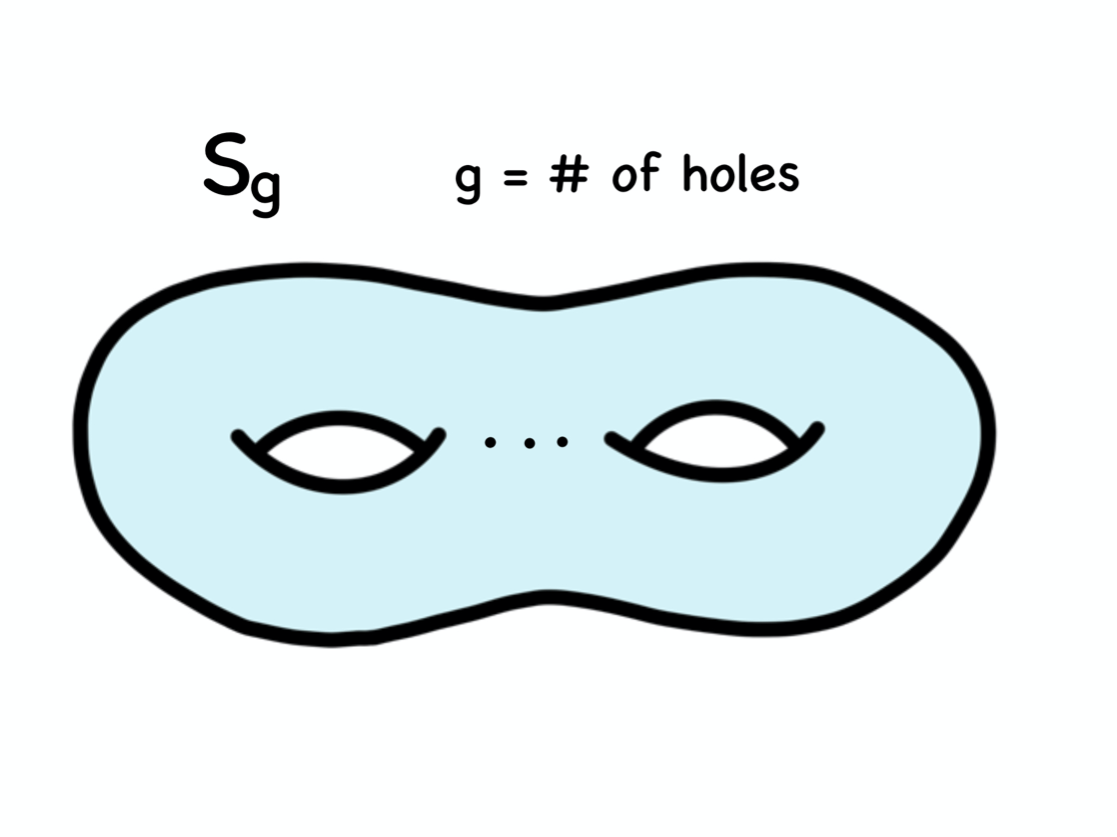}\end{center} 
  \caption{Surface with genus $g\ge 2$}
 \end{figure}
Now that our traveller learned that there exist many (even a continuous family of) closed hyperbolic surfaces to explore, she has to understand what a Euclidean line is in this world.

\section*{Euclidean lines in closed hyperbolic surfaces}
In the upper half plane model of $\bH^2$,
the horizontal line
$\{y=1\}$ will be called a Euclidean line in $\bH^2$.
In a given geometric space, objects which are isometric to each other should be given the same name. We note that the image of $\{y=1\}$ under an isometry of $\bH^2$ is either another horizontal line, or a circle which is tangent to the $x$-axis. So all of these objects will be called Euclidean lines in $\bH^2$. A nickname for a Euclidean line is a horocycle. A horocycle in $\bH^2$ is characterized as an isometric embedding of the real line $\br$ to $\bH^2$ with constant curvature one; geodesics in $\bH^2$ have constant curvature zero.
\begin{figure}[ht] \begin{center}
 \includegraphics [height=3cm]{horo1}\end{center}
 \caption{Euclidean lines in $\bH^2$}\label{horo1}
\end{figure}

From the point of view of our Euclidean traveller, imagine that she wants to drive in a car where the steering wheel is in one fixed position and 
to travel without bumps. If her car is turning at a constant rate, she can lean back against the seat, feeling stable. It's the change in curvature that makes the car trip bumpy. If the wheel is fixed at a small angle, she stays within a bounded distance of a geodesic.
If the wheel is fixed at a large angle, she goes around in a circle and stays a bounded distance from a point.  But in between there is a perfect angle where neither happens, and she moves along a horocycle!  

A Euclidean line in a closed hyperbolic surface
$\Gamma\ba \bH^2$ is the image of a Euclidean line in $\bH^2$ under the quotient map $$\pi: \bH^2\to \Gamma\ba \bH^2.$$

Where does our Euclidean traveller get to visit in a closed hyperbolic surface?

 \begin{figure}[ht] 
  \includegraphics [height=4cm]{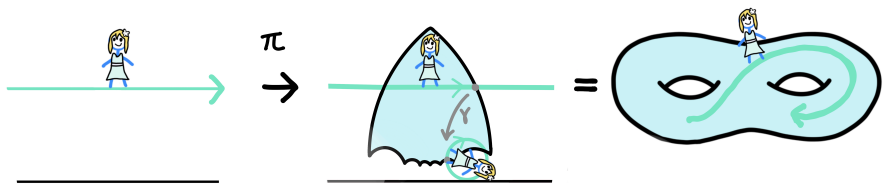} 
\caption{Travelling along a Euclidean line} \label{figT}\end{figure}

Here is an illustration given in Figure \ref{figT}.
Fix a (blue colored) fundamental domain $\mathsf O$ for $\Gamma\ba \bH^2$.
When the traveller, walking along the  mint-colored Euclidean
line $L$, reaches the boundary of
$\mathsf O$, she  gets instantly moved to a different side of $\mathsf O$ by some side-pairing transformation
$\gamma\in \Gamma$. She then continues her journey on the 
new Euclidean line $\gamma (L)$ inside  $\mathsf O$ until she reaches the boundary of $\mathsf O$ again, etc.
The instant jumps and the shapes of translates of $L$ made by $\Gamma$ in $\bH^2$ appear more complicated than those in the Euclidean torus $\mathbb T^2$.

Nevertheless, Hedlund (1936) assures that our Euclidean traveller gets to see all the places in a closed hyperbolic surface, no matter where her initial point of departure is.

\begin{thm} \cite{He}\label{He}
 Any Euclidean line in a closed hyperbolic surface is dense.
\end{thm}

\begin{figure}[ht] \begin{center}
\includegraphics [height=4cm]{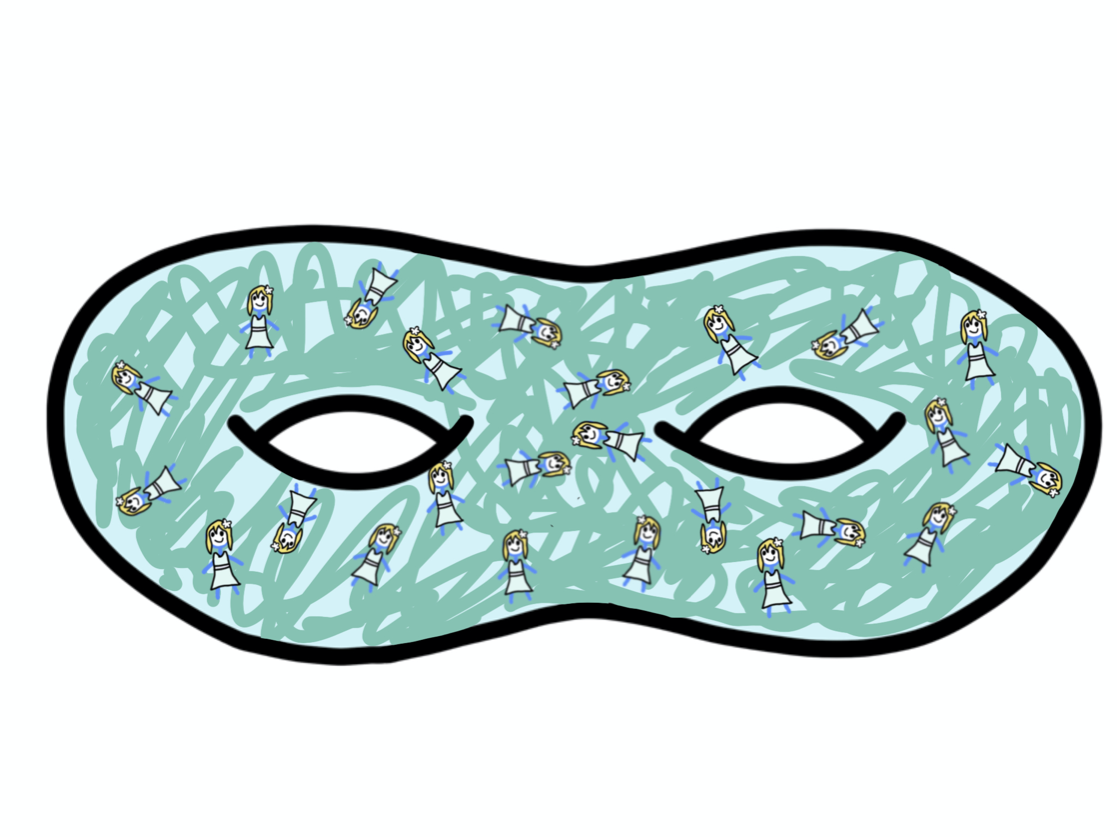}\end{center}
  \caption{Euclidean traveller} \label{den}
\end{figure} 

\subsection*{Hyperbolic lines can be very wild}
We remark that this theorem of Hedlund is about {\it Euclidean lines}.
The closure of a hyperbolic line does not even have to be a submanifold in general.
 \begin{figure}[ht] \begin{center}
  \includegraphics [height=4cm]{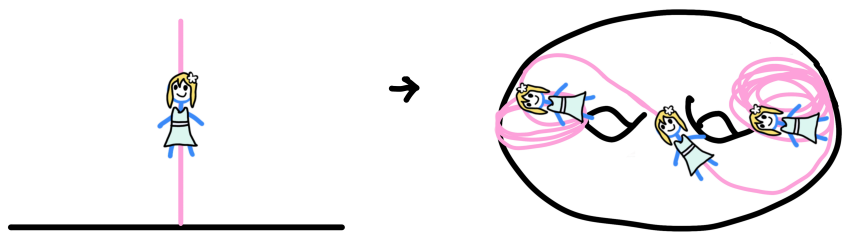} 
  \caption{Hyperbolic traveller}\label{HT} \end{center}
\end{figure}
 As illustrated in Figure \ref{HT}, a geodesic can ``spiral'' around closed geodesics and its closure can be a fractal of {\it dimension} strictly between $1$ and $2$.
\subsection*{Going upstairs}
\begin{figure}[ht] \begin{center}
  \includegraphics [height=3cm]{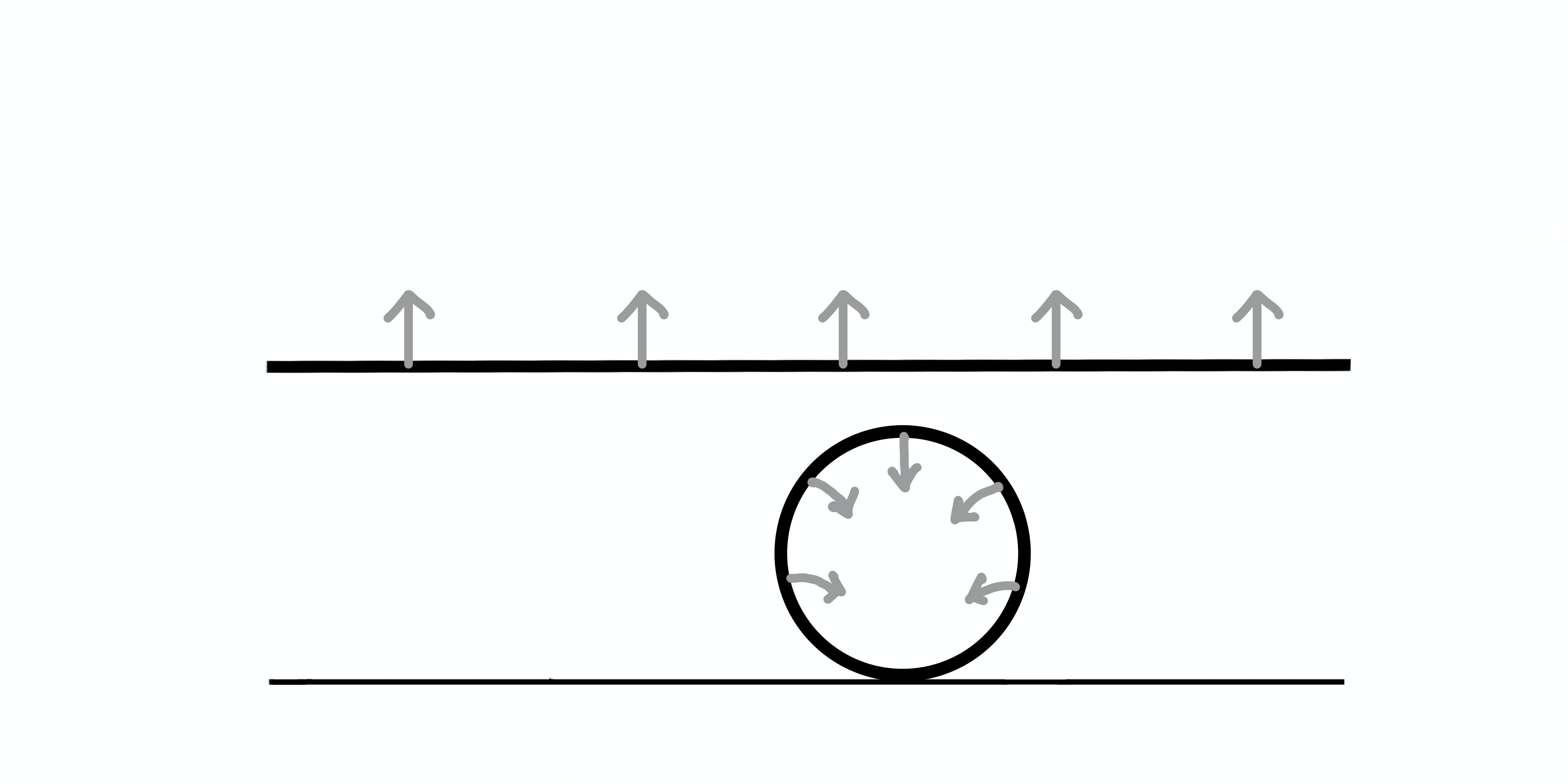}    
  \end{center} \caption{$U$-orbits}\label{ho}
\end{figure}

Hedlund's proof of Theorem \ref{He} relies on the fact that
 $\bH^2$ is almost $\PSL_2(\br)$. The isometric action of $\PSL_2(\br)$ on $\bH^2$,
 which gave us the identification $\bH^2=\PSL_2(\br)/\op{SO}(2)$, extends to an action
 on the unit tangent bundle $\T^1(\bH^2)$. In this action, the stabilizer of a vector is trivial, as no rotation in the plane fixes  a vector. If we denote $v_0\in \T^1(\bH^2)$ the upward normal vector based at $i$, then the orbit map $g\mapsto g(v_0)$ now gives an isomorphism $\PSL_2(\br)=\T^1(\bH^2)$, where the identity matrix of $\PSL_2(\br)$
 corresponds to the vector $v_0$.
 Moreover, if we consider the following one dimensional subgroup
$$U=\left\{\begin{pmatrix} 1& t \\ 0 &1\end{pmatrix}: t\in \br\right\},$$
 then the subgroup $U$ corresponds to
 the set of all upward normal vectors on the Euclidean line $\{y=1\}$ (see Figure \ref{ho}).

Similarly, for $g\in \PSL_2(\br)$ with $g\{y=1\}$ a circle tangent at $\xi$, 
the  coset $gU$ corresponds to the set of all inward  normal vectors on $g\{y=1\}$ pointing to  $\xi$.  Indeed, any Euclidean line (EL) in $\bH^2$ arises as the image of some $gU$ under the basepoint projection map $\PSL_2(\br)=\op{T}^1(\bH^2)\to \bH^2$ given by $g\mapsto g(i)$:

\begin{equation*}\boxed{
\begin{array}{ccccc}
   \op{T}^1( \bH^2) & = & \PSL_2(\br)  & \supset  gU  \\
   & \downarrow  &  & \downarrow   \\ 
 \bH^2  &=&  \PSL_2(\br)/\op{SO}(2)  & \supset  \text{EL} \\
\end{array} } \end{equation*}

Since the image of $gU$ in $\Gamma\ba \PSL_2(\br)$, under the projection $\PSL_2(\br)\to \Gamma\ba \PSL_2(\br)$, is equal to $\Gamma\ba \Gamma gU$
and $\bH^2=\PSL_2(\br)/\op{SO}(2)$,
this picture is preserved under the quotient map $\pi:\bH^2\to \Gamma\ba \bH^2$:

\begin{equation*}\boxed{
\begin{array}{cccc}
   \op{T}^1(\Ga\ba  \bH^2) & = & \Ga\ba \PSL_2(\br)  & \supset  xU  \\
   & \downarrow  &  & \downarrow   \\ 
 \Ga\ba \bH^2  &=&  \Ga\ba \PSL_2(\br)/\op{SO}(2)  & \supset  \text{EL} \\
\end{array} } \end{equation*}

It follows that any Euclidean line in $\Gamma\ba \bH^2$ 
is of the form $xU (i)$ and $\overline{xU(i)}= \overline{xU} (i)$, as the basepoint projection map $g\to g(i)$ has compact fibers. Therefore
 if we can describe the closures of all orbits $xU$ in $\Gamma\ba G$,
we understand the closure of a Euclidean line in $\Gamma\ba \bH^2$.

\begin{figure}[ht] \begin{center}
  \includegraphics [height=4cm]{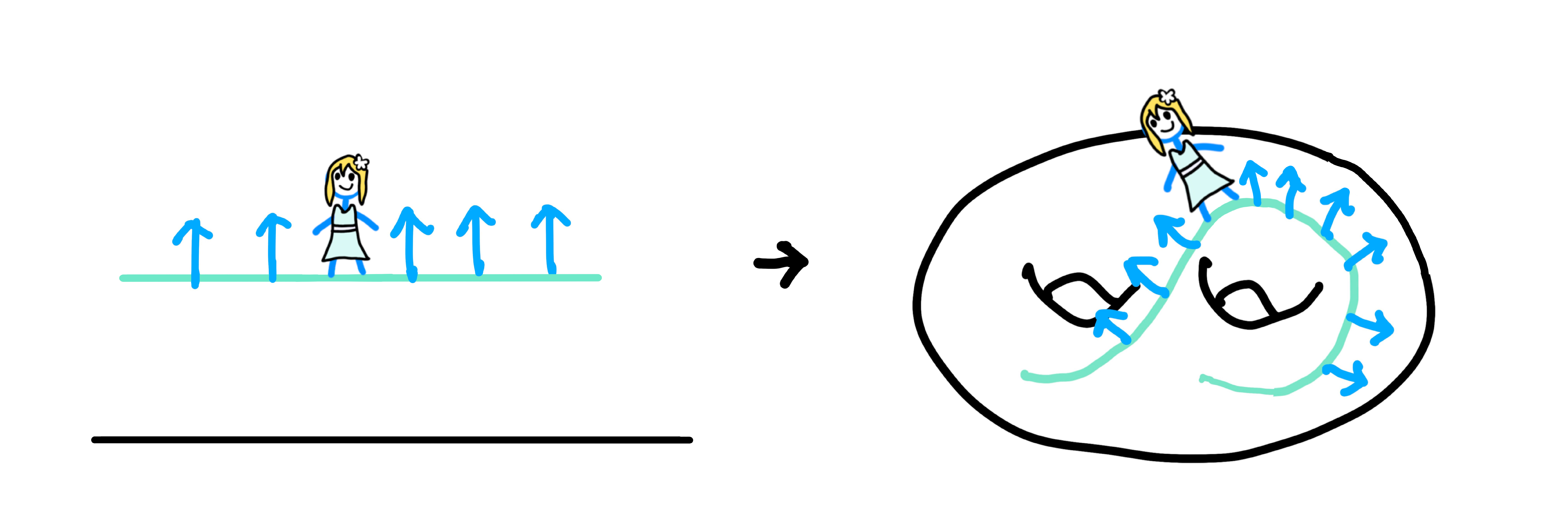} 
  \end{center}\caption{Travelling with an arrow}
\end{figure}

Indeed, Hedlund says that every $U$-orbit is dense upstairs in $\op{T}^1(\Gamma\ba \bH^2)$:
\begin{thm} \cite{He} For any $x\in \Ga\ba \PSL_2(\br)$,
$$\overline{xU}=\Ga\ba \PSL_2(\br).$$
\end{thm}

This means that our Euclidean traveller in a closed hyperbolic surface not only gets to see all the places, but she is able to appreciate those places from all angles as in Figure \ref{den}. 

\section*{Closed hyperbolic $n$-manifolds}
For $n\ge 2$, the hyperbolic $n$-space $\bH^n$ is the unique simply connected complete  $n$-dimensional Riemannian manifold of
constant sectional curvature $-1$. Its upper half space model is given by
$$\bH^n=\{(x_1, \cdots, x_{n-1},y): y>0\}$$
with the hyperbolic metric $ds=\tfrac{\sqrt{dx_1^2+\cdots +dx_{n-1}^2+dy^2}}{y} .$ In this model, geodesics are  vertical lines or semicircles meeting the hyperplane $\br^{n-1}\times \{0\} $ orthogonally. 

The hyperboloid model (also called the Minkowski model) of $\bH^n$ is given by
$\{Q_0(x_1, \cdots, x_{n+1})=-1, x_{n+1}>0\} $
where 
\begin{equation}\label{q} Q_0(x_1, \cdots, x_{n+1})=x_1^2 +x_2^2+\cdots + x_n^2 -x_{n+1}^2 .\end{equation}
In this model, the hyperbolic distance $d(x, y)$ is given by
$\cosh d(x,y)= -Q_0(x,y)$ (cf. \cite{Ratc}). In particular, any element $g$ of 
the special orthogonal group
$\op{SO}(Q_0)$, that is,
$g\in \op{SL}_{n+1}(\br)$ satisfying
$$Q_0(g v)= Q_0(v)\text{ for all $v\in \br^{n+1}$},$$
is an isometry of $\bH^n$. Indeed,
the group $\op{Isom}^+(\bH^n)$
of orientation preserving isometries is given by the identity component of the special orthogonal group $\op{SO}(Q_0)$; so
$$\op{Isom}^+(\bH^n)\simeq \op{SO}^\circ (n,1).$$

This is consistent with our previous statement that
$\op{Isom}^+(\bH^2)=\op{PSL}_2(\br)$, since $\op{PSL_2}(\br)\simeq \op{SO}^\circ (2,1)$.
As in the dimension $2$ case, we have
\bi
\item any closed hyperbolic $n$-dimensional manifold is
a quotient $$M=\Ga\ba \bH^n$$ where $\Ga$ is a discrete (torsion-free) cocompact subgroup
of $\op{SO}^\circ (n,1)$.
\ei

Unlike the dimension two case where there is a continuous family of closed hyperbolic surfaces, higher dimensional closed hyperbolic manifolds are rarer.
The Mostow rigidity theorem \cite{Mos} implies that there exist
only countably many closed hyperbolic manifolds of dimension at least $3$.
Nevertheless, there are infinitely many of them \cite{Bo}.

In the upper half-space model $\bH^n=\{(x_1, \cdots, x_{n-1},y): y>0\}$,
the horizontal line $\{(x_1, 0, \cdots, 0, 1): x_1\in \br\}$, and its isometric images
will be called a Euclidean line (=horocycle) in $\bH^n$. As before, a Euclidean line in $M=\Gamma\ba \bH^n$ is 
the image of an Euclidean line in $\bH^n$ under the canonical quotient map $\pi: \bH^n\to \Gamma\ba \bH^n$.  They are isometric immersions of $\br$
with zero torsion and constant curvature $1$.

\bigskip 
{\it{What places does our Euclidean traveller get to see
in $M$?}}

\begin{figure}[ht] \begin{center}
  \includegraphics [height=4cm]{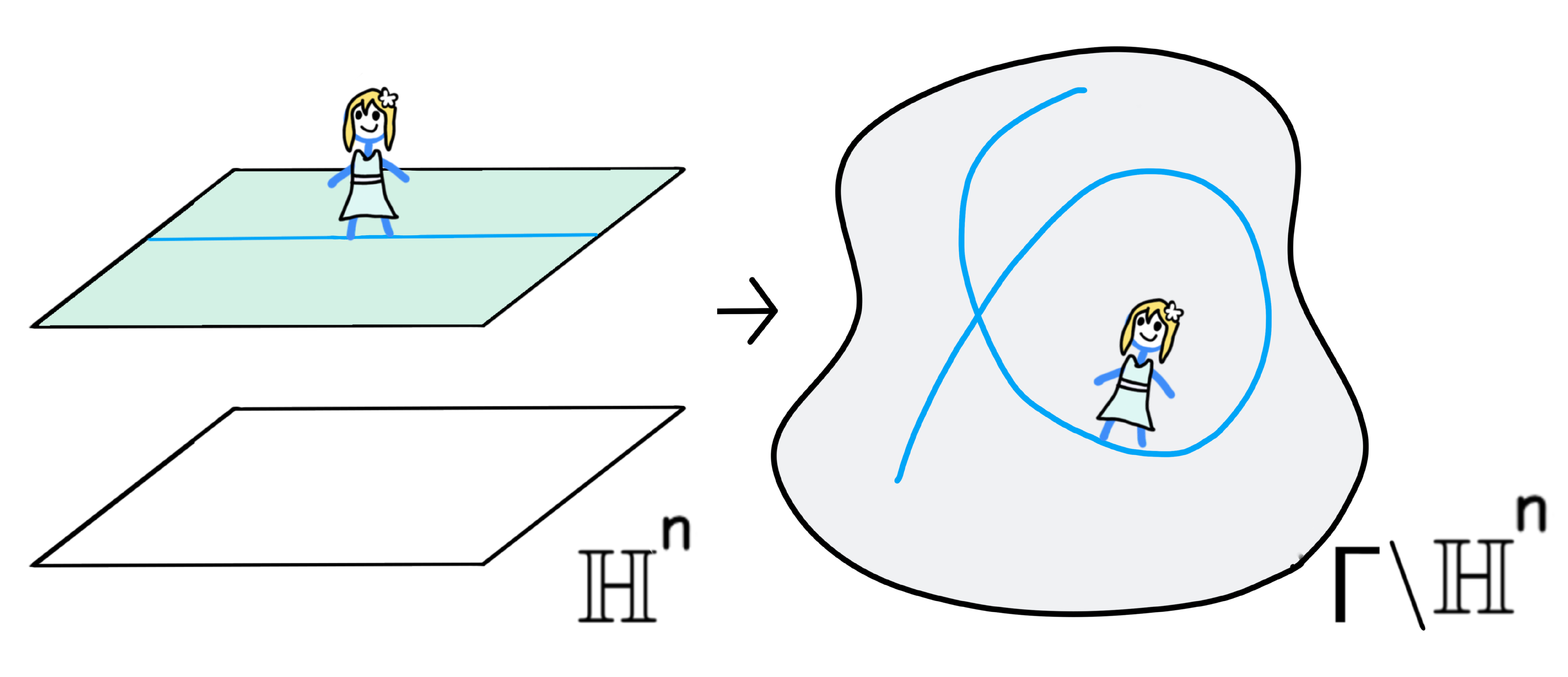}
  \end{center}
  \caption{Euclidean traveller in $\Gamma\ba \bH^n$}
\end{figure}

The answer to this question is a special case of Ratner's theorem on orbit closure classification.

\section*{Homogeneous dynamics}
Let $G$ be a connected simple linear Lie group.
(e.g., $\op{SL}_n(\br), \op{SO}^\circ(n,1)$), and  
$\Ga<G$ be a discrete subgroup.
The quotient space $\Ga\ba G$ is a homogeneous space, in the sense that there
is a transitive action of a Lie group, which is $G$ in our case. 
Any subgroup $U$ of $G$ acts on this homogeneous space
$\Ga\ba G$ by translations on the right, giving rise to a topological dynamical system:
$$\Gamma\ba G \curvearrowleft U.$$
Studying dynamical properties of subgroup action on a homogeneous space
$\Gamma\ba G$ is a main subject in the field of {\it homogeneous dynamics}.
A fundamental question in homogeneous dynamics is whether one can understand all possible orbit closures:
\begin{Q}
For a given point $x\in \Ga\ba G$, what is the closure $\overline{xU}?$
\end{Q}
{} 

 Moore's ergodicity theorem \cite{Mo} implies that
if the homogeneous space $\Gamma\ba G$ is compact, or more generally is
of  finite volume, then for any non-compact subgroup
 $U$ of $G$, almost all $U$-orbits are dense, that is, for almost all $x\in \Ga\ba G$,
 $$\overline{xU}=\Ga\ba G .$$ 
 Indeed, a hyperbolic line in a closed hyperbolic surface
 $\Gamma\ba \bH^2$ is the image of an orbit of the diagonal subgroup
 $A=\{\tiny{\begin{pmatrix} e^{t} & 0\\ 0&e^{-t}\end{pmatrix}}:t\in \br \}$ of $\PSL_2(\br)$ under the basepoint projection $\Gamma\ba \PSL_2(\br)=\T^1(\Gamma\ba \bH^2)\to \Gamma\ba \bH^2$; hence
 this theorem implies that almost all  hyperbolic lines are dense in a closed hyperbolic surface.
 
The weakness of this ergodicity theorem is that it is only the existence of many dense orbits. If one specifies the initial point $x$,  it gives no information about  the closure of the given orbit $xU$.

The following celebrated theorem of Ratner fixes this weakness for {\it unipotent} subgroup actions, answering the conjecture of
Raghunathan in the affirmative.

A matrix $u$ is called {\it unipotent} if all of its eigenvalues are equal to $1$. For instance, $u_t=\begin{pmatrix} 1& t \\ 0 &1\end{pmatrix} $ is unipotent
but $a_t=\begin{pmatrix} e^{t} & 0\\ 0&e^{-t}\end{pmatrix}$ is not unipotent for $t\ne 0$.

\begin{thm}\cite{Ra} \label{Ra} Let $\op{Vol}(\G\ba G)<\infty$ and $U$ be a connected subgroup of $G$ generated by unipotent elements.
For any $x\in \Gamma\ba G$,
$$\overline{xU}=xH $$
for some closed connected subgroup $H$ of $G$ which contains $U$.
\end{thm}

\section*{Back to Euclidean lines in closed hyperbolic manifolds}
Given a closed  hyperbolic manifold $M=\Gamma\ba \bH^n$, its frame bundle $\op{F}(M)$, consisting of oriented orthonomal frames on $M$, is identified with the homogeneous space $\Gamma\ba \op{SO}^\circ (n,1)$, as the action of $\op{Isom}^+(\bH^n)=\op{SO}^\circ (n,1)$ is simply transitive on $\op{F}(\bH^n)$.
Moreover, we have a one-dimensional subgroup $U$ of $\op{SO}^\circ (n,1)$ given by
   $$U=\left\{\begin{pmatrix}  1& t & 0 & \cdots & 0 & -t^2/2\\
   0 & 1 & 0 &\cdots & 0 & -t\\
   0 & 0 & 1 &\cdots & 0 & 0 \\
   \vdots & \vdots & \vdots & \vdots & \vdots & \vdots \\
   0 &0 &0 &0 &1 & 0\\
   0 &0 &0 &0 &0 & 1\end{pmatrix} :t\in \br \right\}$$
   (up to conjugation) such that
    every Euclidean line (EL) in $M$ is the image of some $U$-orbit in
    $\Gamma\ba \op{SO}^\circ (n,1)$ under the basepoint projection map:
    
    \begin{equation*}\boxed{
\begin{array}{ccccc}
 \op{F}(\Ga\ba \bH^n) & = & \Gamma\ba \op{SO}^\circ (n,1)  & \supset  xU  \\
   & \downarrow  &  & \downarrow   \\ 
  & \Gamma\ba \bH^n&   & \supset  \text{EL} \\
\end{array} } \end{equation*}

Since $U$ consists of unipotent elements, Ratner's theorem \ref{Ra} 
applies to orbits of $U$.

For each integer $1\le k\le n$, we set $\mathcal H_k:=\{(x_1, \cdots, x_{k-1}, 0, \cdots, 0, y)\in \bH^n\}$.  A hyperbolic $k$-subspace of $\bH^n$, which will be denoted by $\bH^k$, is an isometric image of $\mathcal H_k$, that is, $g (\mathcal H_k)$ for some $g\in \op{Isom}^+(\bH^n)$. 

\begin{Def} By a closed hyperbolic (properly immersed) $k$-submanifold of $M$, we mean the
compact image 
$$\pi(\bH^k)\subset \Gamma\ba \bH^n$$ for a hyperbolic $k$-subspace $\bH^k\subset \bH^n$. \end{Def}
Note that for any $k\ge 2$, a closed hyperbolic $k$-submanifold of $M$ contains many dense Euclidean lines.
So certainly, if $M$ possesses a closed hyperbolic $k$-submanifold $N$, then some Euclidean lines lying in $N$ will have their closures equal to $N$. There are other possibilities for the closure of a Euclidean line in $M$. To explain what they are, we introduce the notion of a {\it tilting} of a hyperbolic subspace, which is analogous to a translation of a subspace
by a vector in the Euclidean space $\br^n$.

\begin{figure}[ht]
 \begin{center} \includegraphics [height=4cm]{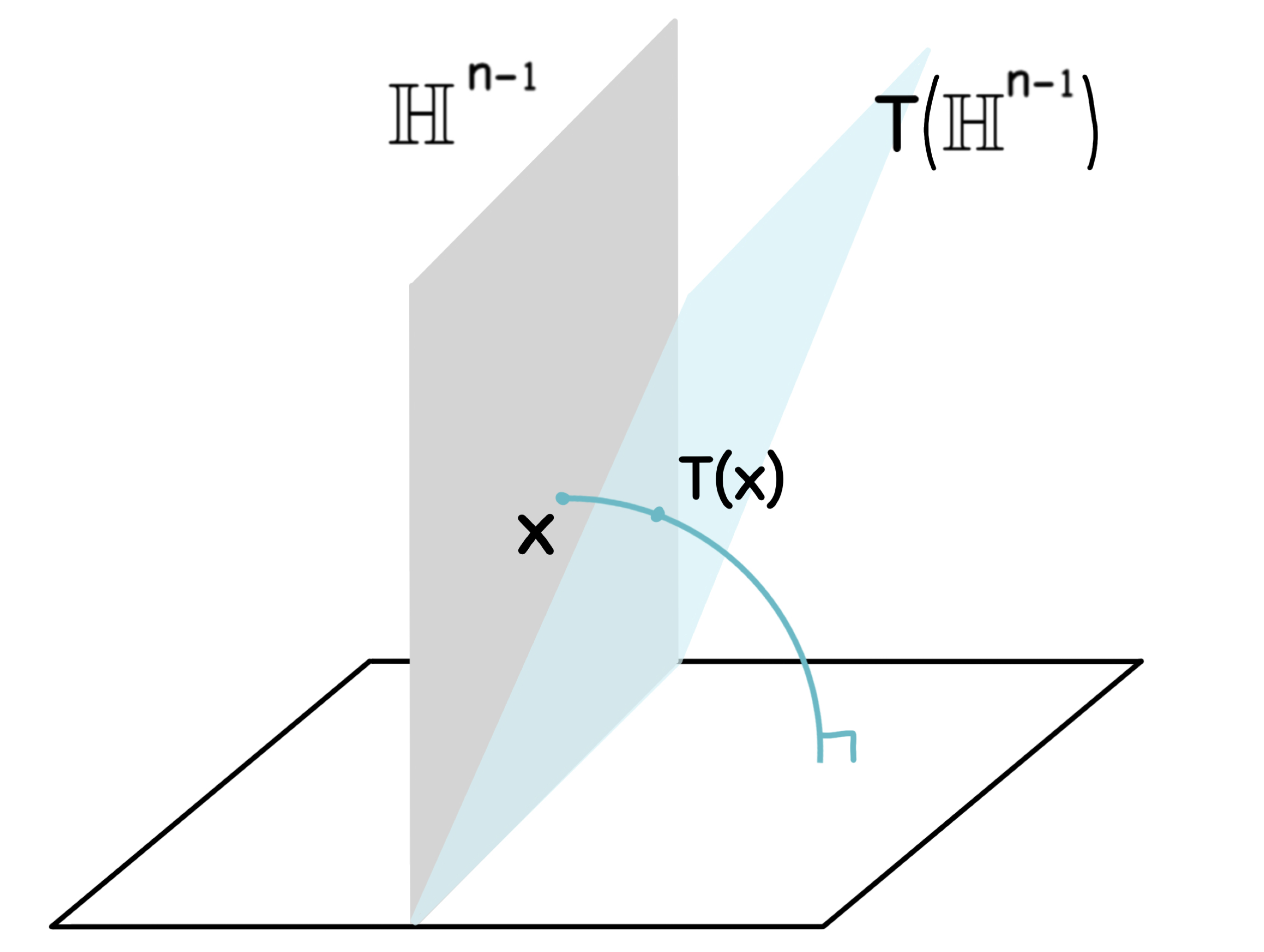} \end{center}
 \caption{Tilted hyperbolic subspace}\label{tilt}
\end{figure}

\begin{Def}[Tilting] \rm 
For a given hyperbolic $k$-subspace $\bH^k$ of $\bH^n$ for $1\le k\le n-1$, choose a hyperbolic subspace $\bH^{k+1}$ containing $\bH^k$ equipped with an orientation, and a nonnegative number $d_0\ge 0$. For $x\in \bH^{k}$,
let $\mathsf{T} (x)\in \bH^{k+1}$ be the unique point lying in the geodesic determined by the outward normal vector to
$\bH^{k}$ and $d(x,\mathsf{T} (x))=d_0$. We call $\mathsf{T}:\bH^k\to \bH^{k+1}$ a tilt map and
its image $\mathsf{T}(\bH^{k})$ a tilted hyperbolic $k$-subspace of $\bH^n$ (Figure \ref{tilt}).
\end{Def}
The significance of a tilt map relevant to our discussion is that the image of a Euclidean line lying in $\bH^k$ under a tilt map is a Euclidean line and that if $N=\pi(\bH^k)$ is a closed hyperbolic submanifold of $M$, then $\pi(\mathsf{ T}(\bH^k))$ is also a compact submanifold
which is equidistant from $\pi(\bH^k)$.
Noting that a tilt map commutes with isometric action of $G$, the following
$$\mathsf{ T} (N)= \pi(\mathsf{ T}(\bH^k))$$ is well-defined;
we call $\mathsf{ T} (N)$ a tilting of a hyperbolic manifold $N$.
If a Euclidean line $L$ is dense in $N=\pi(\mathbb H^k)$, then $\mathsf{ T}(L)$
is dense in $\mathsf{ T} (N)$.

Now the following consequence of Ratner's theorem \cite{Ra} implies that these are all possible closures of a Euclidean line.

\begin{figure}[ht]
 \begin{center}
    \includegraphics[height=3cm]{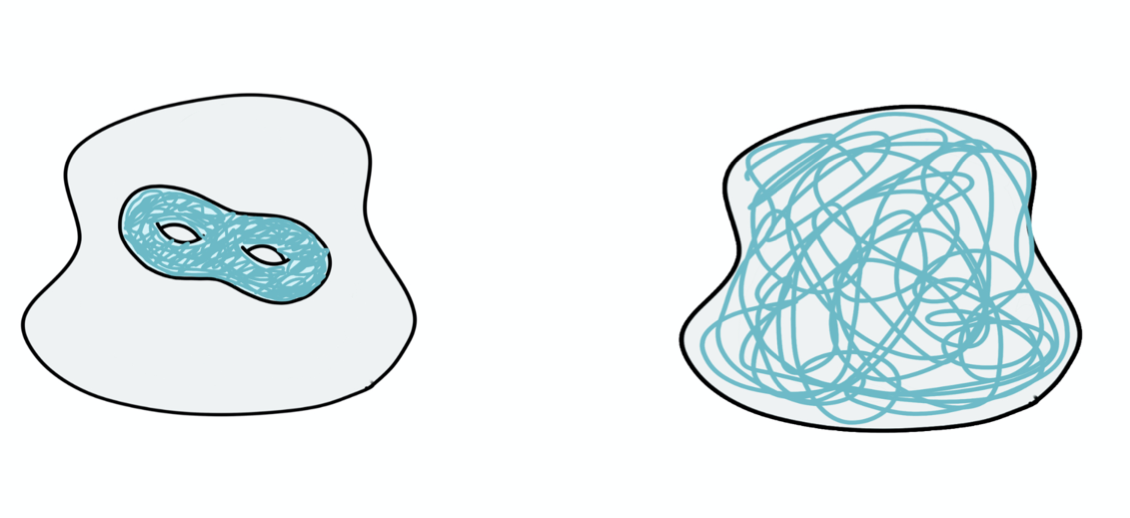} 
 \end{center}\caption{Closures of Euclidean lines}
\end{figure}
\begin{Thm} \label{Rac}   Let $M=\Ga\ba \bH^n$ be a closed hyperbolic
$n$-manifold for $n\ge 2$. Then the closure of any Euclidean line in $M$
is a closed hyperbolic submanifold, up to tilting.
\end{Thm}

Therefore, no matter where she starts her journey, our Euclidean traveller in a closed hyperbolic world is guaranteed to see all the places of some  hyperbolic subworld
or at least a tilted version of it.

\subsection*{Travelling only in one direction} 
From a traveller's point of view, a natural question to the travel guide is whether she can still sightsee the same places by walking only in one direction,
in other words, is the closure of a half Euclidean line in a closed hyperbolic manifold same as the closure of the whole Euclidean line?
The answer is {\it yes} by the following version of Ratner's theorem:
\begin{thm} \cite{Ra} Let $\op{Vol}{(\Gamma\ba G)}<\infty$ as in Theorem \ref{Ra}, and let $U=\{u_t:t\in \br \}$ be a one-parameter subgroup  of unipotent elements.  For any $x\in \Gamma\ba G$,
\begin{equation}\label{plus} \overline{xU}=\overline{xU^{+}} \end{equation}
where $U^{+}=\{u_{t}: t\ge 0\}$.
\end{thm}

\section*{Hyperbolic manifolds of infinite volume}
After adventuring all closed hyperbolic manifolds, our Euclidean traveller is now ready to venture into a hyperbolic manifold of infinite volume.
Will she again be able to see all the places of some hyperbolic subworld?

It turns out that there are
certain hyperbolic manifolds homeomorphic to the product of a closed surface and $\br$ where
some Euclidean lines have wild closures. On the other hand, our wonderful travel agency presents a list of infinitely many hyperbolic manifolds of infinite volume where similar well-planned sightseeing is possible as in  closed hyperbolic manifolds. This class of
hyperbolic manifolds are called hyperbolic manifolds with {\it Fuchsian ends}.
They are all obtained from closed hyperbolic manifolds by a certain removing and growing process; we may think of them as children of closed hyperbolic manifolds.
\section*{Hyperbolic surfaces with Fuchsian ends}
\begin{figure}[ht]
 \begin{center}
    \includegraphics[height=4cm]{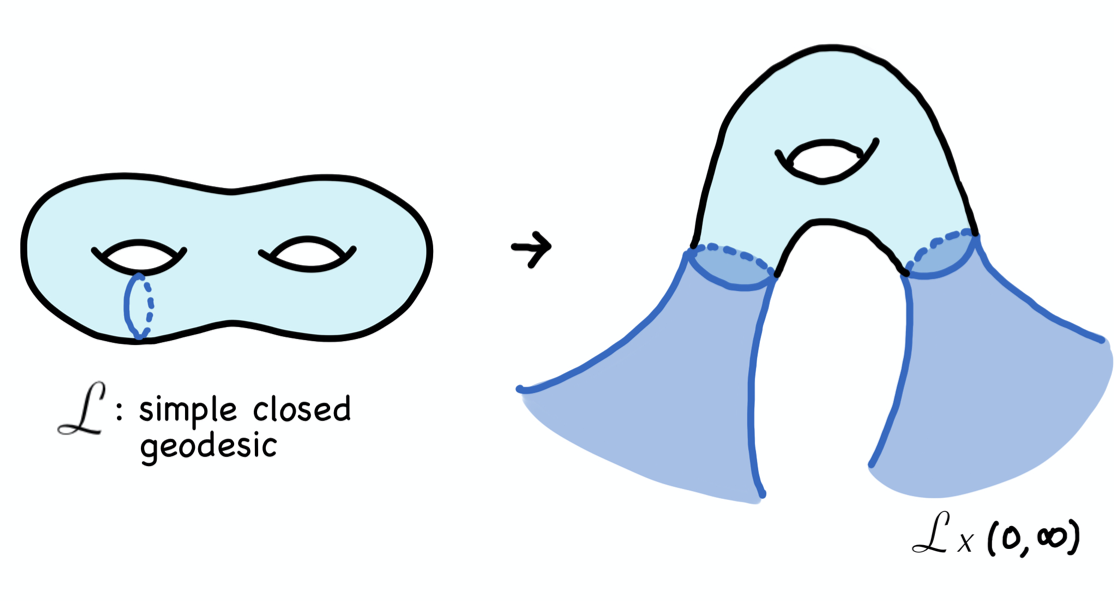} 
 \end{center}
 \caption{Hyperbolic surface with Fuchsian ends}\label{FE0}
\end{figure}

We first explain how to construct a hyperbolic surface with Fuchsian ends as illustrated in Figure \ref{FE0}.
For any closed hyperbolic surface $S=\Gamma_0\ba \bH^2$, choose a simple closed geodesic $\mathcal L$
($S$ contains infinitely many such), which may or may not be separating.

After removing $\mathcal L$ from $S$, we get a hyperbolic surface $N_0$ whose boundary is one or two copies of $\mathcal L$, depending on whether $\mathcal L$ is separating or not. 
If we let each boundary component {\it grow naturally} in hyperbolic world, we get a hyperbolic surface with Fuchsian ends.
In other words, there is a canonical way to extend $N_0$ to a complete hyperbolic surface which is homeomorphic to $N_0$, namely, by gluing  the continuous stack  $\mathcal L\times [0,\infty)=\{\mathcal L^t:
t\in [0, \infty)\}$ of titled hyperbolic lines to each boundary component of $N_0$.
The resulting hyperbolic surface is of the form $N=\Gamma\ba \bH^2$ where $\Gamma<\Gamma_0$ is the fundamental group of $N_0$. This is an example of a hyperbolic surface with Fuchsian ends (Figure \ref{FE0}).

More generally, a hyperbolic surface with Fuchsian ends is obtained from a closed hyperbolic surface $S$ by a similar procedure for a finite set of disjoint simple closed geodesics $\{\mathcal L_1, \cdots, \mathcal L_k\}$.
A connected surface, say $N_0$, of $S-\bigcup_i \mathcal L_i$ is a hyperbolic surface with boundary components and the resulting complete hyperbolic surface, say, $N$, obtained by gluing the corresponding hyperbolic cylinders to each boundary component is a hyperbolic surface with Fuchsian ends.  The metric closure of $N_0$ is a compact hyperbolic surface which is homotopy-equivalent to $N$.
By a hyperbolic surface with Fuchsian ends, we mean a surface obtained in this way.

In any hyperbolic surface with (non-empty) Fuchsian ends, there exist many
Euclidean lines contained in the end of $N$, that is, in $N-N_0$, and they, as well as 
Euclidean lines which are equidistant from them, 
are proper immersions of the Euclidean line in $\bH^2$ via $\pi$.

The following theorem of Dal'bo, which extends Hedlund's theorem \ref{He}, says that they are the only possible non-dense Euclidean lines.
\begin{thm} \cite{Da} If $N$ is a hyperbolic surface 
with Fuchsian ends, any Euclidean line in $N$ is  closed or dense.
\end{thm}
So our Euclidean traveller will simply disappear from the hyperbolic world
or she will be seeing all the places even including the Fuchsian ends.
There is a friendly warning in the pamphlet that
 walking along a half line won't take her to
all the places she would be seeing on the full line unlike closed hyperbolic manifolds.

\section*{Hyperbolic $n$-manifolds with Fuchsian ends} \begin{figure}[ht]
 \begin{center}
    \includegraphics[height=3.5cm]{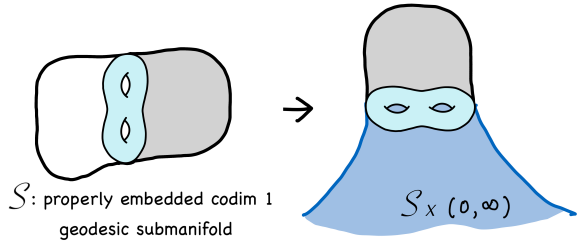}
 \end{center}\caption{Hyperbolic manifold with one Fuchsian end}\label{FE1}
\end{figure}
A hyperbolic $n$-manifold with Fuchsian ends for $n\ge 3$ is constructed similarly
to the surface case, based on the observation that
simple closed geodesics describe all properly embedded hyperbolic submanifolds of a closed hyperbolic surface of codimension one.

Recall that there are only countably many closed hyperbolic $n$-manifolds for $n\ge 3$. 
There are infinitely many of them (though not all) which contain properly embedded hyperbolic submanifolds of codimension one. Take one such closed hyperbolic $n$-manifold $M$
with a properly embedded hyperbolic submanifold $\mathcal S$ of codimension one; see the illustration in Figure \ref{FE1}.

A connected component $N_0$ of $M-\mathcal S$ is a hyperbolic manifold with one or two totally geodesic boundary components isometric to $S$.  If we let $N_0$ naturally grow in
the hyperbolic world, or if we glue  the continuous stack $\mathcal S\times [0, \infty)=\{\mathcal S^t: t\in [0, \infty)\}$ of tilted hyperbolic submanifolds to each boundary component of $N_0$, we obtain a
complete hyperbolic $n$-manifold $N$, homeomorphic to its submanifold $N_0$. This is a hyperbolic manifold with Fuchsian ends;
connected components of $N-N_0$ are called Fuchsian ends of $N$.

\begin{figure}[ht]
 \begin{center}
    \includegraphics[height=5cm]{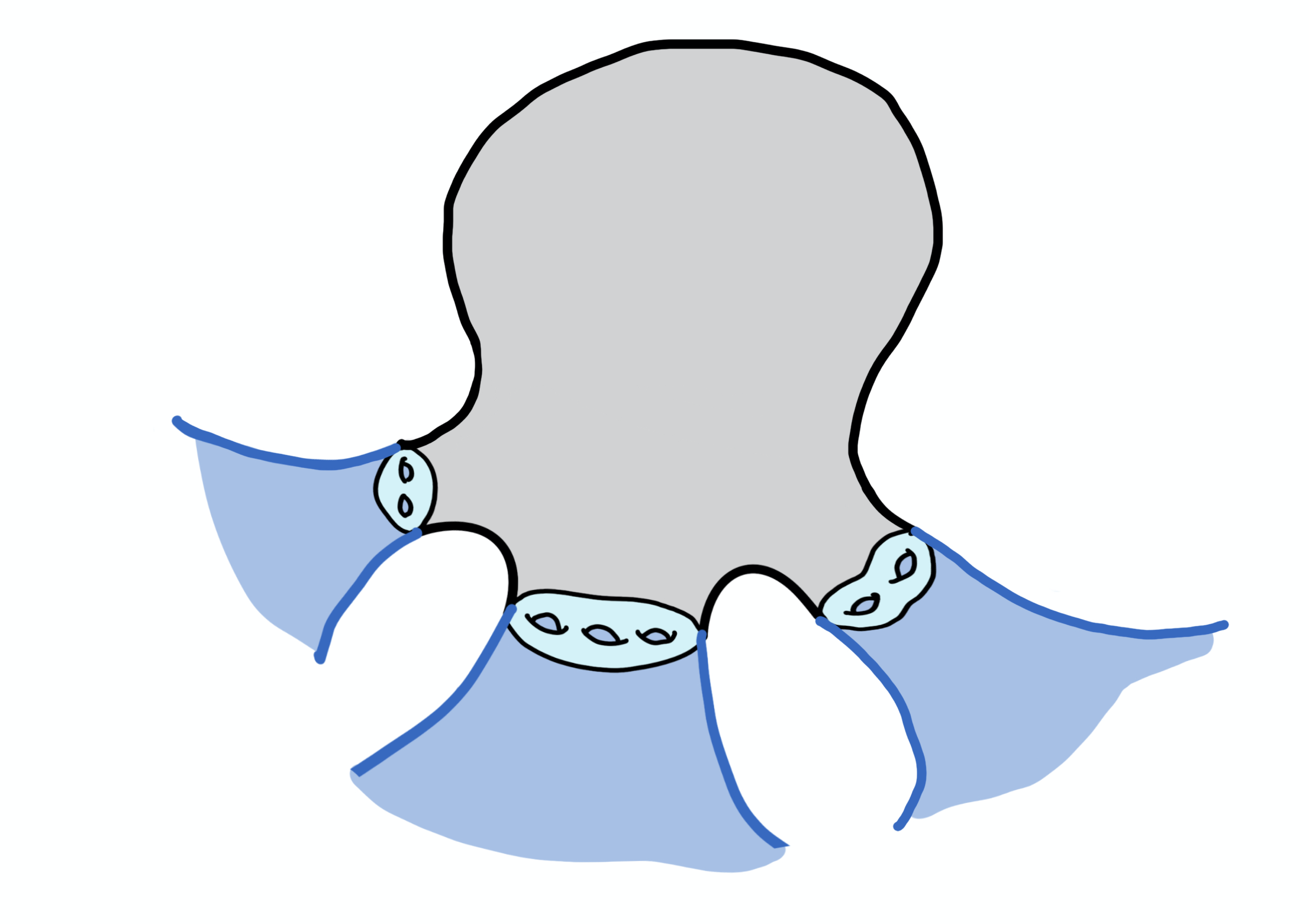}
 \end{center}\caption{Hyperbolic manifold with 3 Fuchsian ends}
\end{figure}

As in the surface case, we mean by a hyperbolic manifold with Fuchsian ends a complete hyperbolic manifold obtained from a closed hyperbolic manifold with a choice of finitely many disjoint properly embedded hyperbolic submanifolds of codimension one.

We may regard a closed hyperbolic manifold as a hyperbolic manifold with empty Fuchsian ends. 
\begin{thm} [\cite{MMO1}, \cite{MMO2}, \cite{LO}]
 Let $N=\Ga\ba \bH^n$ be a hyperbolic manifold with Fuchsian ends for $n\ge 2$. Any Euclidean line in $N$ is closed or its closure is
 a hyperbolic submanifold with Fuchsian ends, up to tilting.
\end{thm}

So our Euclidean traveller in a hyperbolic manifold with Fuchsian ends again finds herself either disappearing from the hyperbolic world or enjoys her sightseeing in some hyperbolic submanifold with Fuchsian ends, or at least a tilted version of it. 

\bigskip 
\noindent{\bf{Acknowledgements}} This article is based on the {\it Erd\"os lecture for Students} that the author gave  at the Joint Mathematics Meetings, 2022. She would like to thank Joy Kim for her help with pictures and Mikey Chow, Dongryul M. Kim, Curt McMullen and Yair Minsky for useful comments on the preliminary version.

\bibliographystyle{plain} 
\bibliography{ET}

       \bigskip
       \bigskip
\begin{tabular}{cl}  
         \begin{tabular}{c}
           \includegraphics[height=3cm]{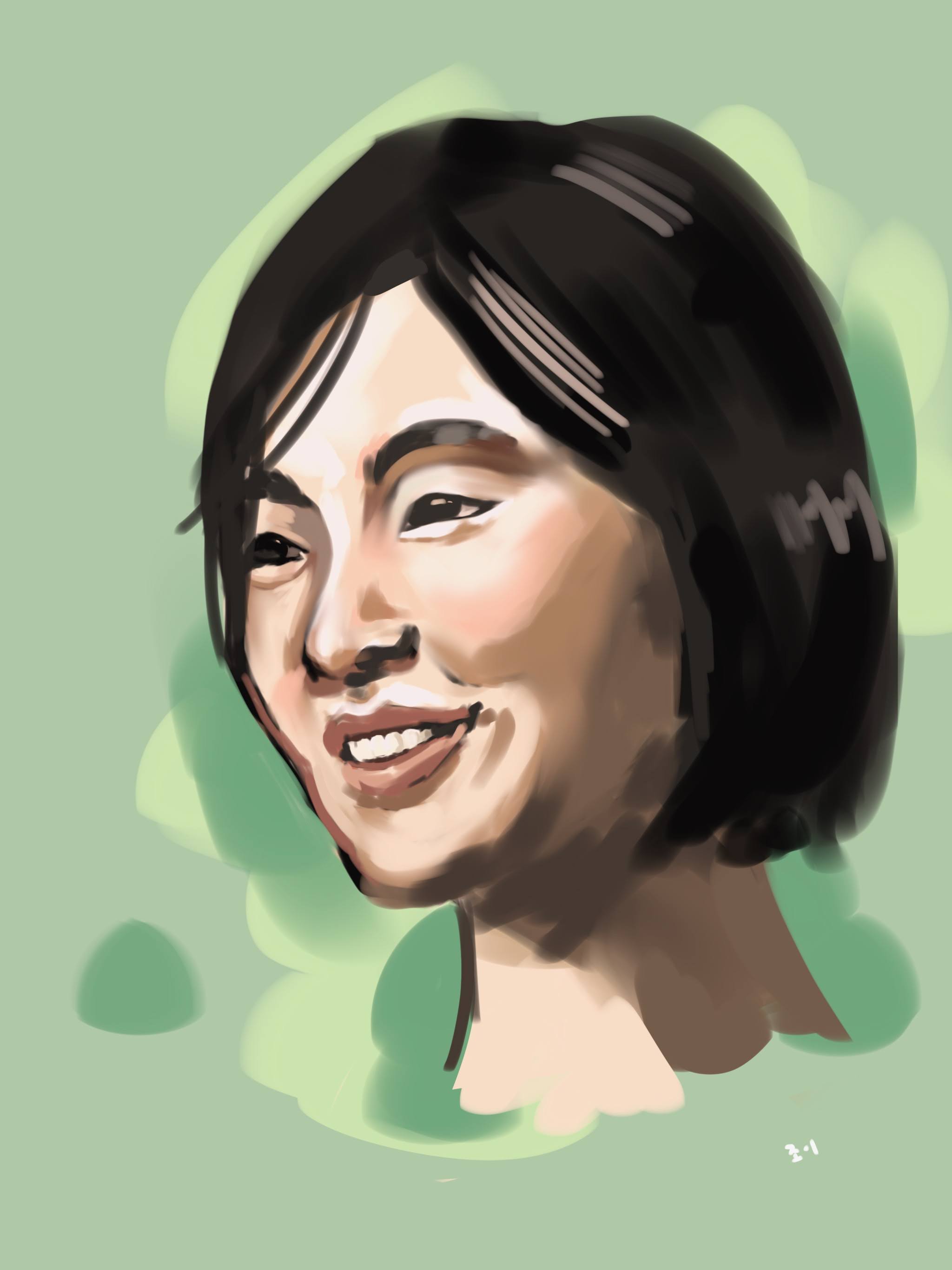}
           \end{tabular}
           & \begin{tabular}{r}
             \parbox{0.5\linewidth}{Hee Oh is the Abraham Robinson Professor of Mathematics at Yale University. Her email address is hee.oh@yale.edu. Picture credit: Joy Kim.}
         \end{tabular} \\
\end{tabular}

\end{document}